%% file: cm.tex
\documentstyle[12pt]{article}
\begin{document}
\title{
The Free Boundary Problem in the \\
 Optimization of Composite Membranes }
\author{S. Chanillo, D. Grieser,  K. Kurata}
\date{\today}
\maketitle
\begin{abstract}
In this paper, continuing our earlier article [CGIKO], 
 we study qualitative properties of 
solutions of a certain eigenvalue optimization problem. 
Especially we focus on  the study of 
the free boundary of our optimal solutions on general domains. 
\end{abstract}
\newtheorem{th}{Theorem}
\newtheorem{prop}{Proposition}
\newtheorem{cor}{Corollary}
\newtheorem{lem}{Lemma}
\newtheorem{prob}{Problem}

\input cm1 
\input cm2 
\input cm3 

\input cm4  
\input cm5 

\vspace{5mm}

\input cm-ref 

\vspace{5mm}

{\sc Address}:

Sagun Chanillo

Department of Mathematics

Rutgers University

New Brunswick, NJ 08903, USA

chanillo@math.rutgers.edu

\medskip
Daniel Grieser

Department of Mathematics

Humboldt-Universit\"{a}t Berlin

Unter den Linden 6

10099 Berlin, Germany

grieser@mathematik.hu-berlin.de

\medskip
Kazuhiro Kurata

Department of Mathematics

Tokyo Metropolitan University

Minami-Ohsawa 1-1, Hachioji-shi, Tokyo, Japan

e-mail: kurata@comp.metro-u.ac.jp

\end{document}

%% file: cm1.tex
\section{Introduction and Summary of results}
In this note, we study some qualitative properties 
of solutions of a certain eigenvalue optimization problem. 
Our note is a continuation and also a summarization of the key results 
of our earlier article [CGIKO]. 
Our problem can be stated in physical terms as :

\noindent 
{\bf Problem(P)}\ Build a body of a prescribed shape out of given materials 
of varying densities, in such a way that the body has a prescribed mass and 
with the property that the fundamental frequency of the resulting membrane 
(with fixed boundary) is lowest possible. 

The physical problem can be re-formulated as a more general 
mathematical problem. More precisely, we are given $\Omega\subset {\bf R}^n$, 
a bounded domain with Lipschitz boundaries and numbers $\alpha >0, 
A\in [0, |\Omega|]$ (with $|\cdot |$ denoting volume).  
For any measurable set $D\subset\Omega$, let $\chi_D$ be the 
characteristic function and $\lambda_{\Omega}(\alpha,D)$ the lowest 
eigenvalue $\lambda$ of the problem,
\begin{eqnarray}
-\Delta u+\alpha\chi_D u&=&\lambda u\quad \mbox{on}\quad \Omega 
\label{eq1}
\\
 u&=&0\quad \mbox{on} \quad \partial \Omega. \nonumber
\end{eqnarray}
Define, 
\begin{equation}
\label{eq2}
\Lambda_{\Omega}(\alpha, A)=\inf_{D\subset\Omega, |D|=A} 
\lambda_{\Omega}(\alpha, D).
\end{equation}
Any minimizer for (\ref{eq2}) will be called an optimal 
configuration for the data $(\Omega, \alpha, A)$. 
If $D$ is an optimal configuration and $u=u_{\alpha, D}$ satisfies 
(\ref{eq1}), then $(u_{\alpha, D}, D)$ will be called an optimal pair (or 
solution). The mathematical problem then reads, 

\noindent
{\bf Problem (M)}\ Study existence, uniqueness and qualitative 
properties of optimal pairs.

We will hereon always work under the nomalization 
\begin{equation}
\label{eq3}
\int_{\Omega} u^2 =1, \ u\geq 0.
\end{equation}
Furthermore, changing $D$ by sets of measure zero does not affect 
$\lambda_{\Omega}(\alpha, D)$ or $u$, thus sets $D$ that differ by sets 
of measure zero will be said to be equal. 

A basic tool that we use to analyse our problem is a variational 
characterization of eigenvalues, precisely, 
$$
\lambda_{\Omega}(\alpha, D)=\inf_{u\in H^1_0(\Omega)}R_{\Omega}(u, D), 
\quad 
R_{\Omega}(u, D)=\frac{\int_{\Omega}|\nabla u|^2 
+\alpha\int_{\Omega}\chi_D u^2}{\int_{\Omega}u^2}.
$$
The minimizer $u$ is well known to exist and is an eigenfunction. Thus, 
for $\Lambda_{\Omega}(\alpha, A)$ we have
$$
\Lambda_{\Omega}(\alpha, A)=\inf_{u\in H^1_0(\Omega), |D|=A}
R_{\Omega}(u, D).
$$
The theorem that follows is basic to the questions 
we hope to treat in this paper. The proof of this theorem is to be found 
in [CGIKO]. 
To state our theorem we will need to introduce some notation. 
First, we will consistently use the notation $\{ u=t\}$ for 
 $\{ x\in\Omega; u(x)=t\}$, and $\{ u\le t\}$ for 
$\{ x\in\Omega; u(x)\le t\}$. 
\begin{th}{\sc ([CGIKO])}
\label{th1}
For any $\alpha >0$ and $A\in [0, |\Omega|]$, there exists an optimal pair. 
Moreover any optimal pair $(u_{\alpha, D}, D)$ has the following properties:

\noindent
(a) \ $u_{\alpha, D}\in C^{1, \delta}(\Omega)\cap W^{2,2}(\Omega)\cap 
C^{\gamma}(\overline{\Omega})$ for every $\delta <1$ and some $\gamma >0$.

\noindent
(b)\ $D$ is a sub-level set of $u_{\alpha, D}$, i.e. there exists $t >0$  
such that 
$$
D=\{ x\in \Omega; u_{\alpha, D}(x)\le t\}.
$$

\noindent
(c)\ The value of $\alpha$ for which
$\Lambda_{\Omega}(\alpha,A)=
\alpha,$ is unique.

\noindent
(d)\ Every level set $\{ u_{\alpha, D}=s\}, s\neq t$ has measure zero. 
If in addtion $\Lambda_{\Omega}(\alpha, A)\neq \alpha$, the free boundary, 
${\cal F}=\{ u_{\alpha, D}=t\}$ has measure zero. 
\end{th}
 
  We use the notation $\overline{\alpha}_{\Omega}(A)$ for the unique
value of $\alpha$ in part (c) above, that is,
\begin{equation}
\label{eq4}
\overline{\alpha}_{\Omega}(A)=\Lambda_{\Omega}(\overline{\alpha}_{\Omega}(A),A)
\end{equation}

We also see right away that our problem to determine an optimal pair 
$(u_{\alpha, D}, D)$ 
which seemed linear is now a non-linear problem, 
\begin{eqnarray}
\label{eq5}
-\Delta u_{\alpha, D}+\alpha\chi_{\{ u_{\alpha, D}\le t\}}u_{\alpha, D}
&=&\Lambda_{\Omega}(\alpha, A)u_{\alpha, D} \quad \mbox{on}\quad \Omega\\
u_{\alpha, D}&=&0\quad \mbox{on}\quad \partial \Omega.\nonumber
\end{eqnarray}
Another important remark is that because any 
optimal configuration $D$ is a 
sub-level set and $u_{\alpha, D}=0$ 
on $\partial\Omega$, the set $D$ will always contain 
a tubular neighborhood of $\partial\Omega$, i.e. $D$ always 
contains a boundary layer.

For notational convenience, from now on we will drop the subscript $\Omega$ 
and write  $\Lambda(\alpha, A)$ for $\Lambda_{\Omega}(\alpha, A)$, 
and use the simpler notation $u$ for $u_{\alpha, D}$, in situations 
where no confusion arises. $\|\cdot\|_{\infty}$ will denote the 
supremum norm in $L^{\infty}(\Omega)$ and $\|\cdot\|_2$ the norm in 
$L^2(\Omega)$.

The main focus in this paper will be on the free boundary $\{ u=t\}$ on 
a general domain $\Omega$. 

Before we state the theorems that we prove in this paper, we continue
summarizing some of the salient results of [CGIKO]. 
It is proved there that problem (M) generalizes problem (P). In fact 
for $\alpha\le \overline{\alpha}_{\Omega}(A)$, the solutions of problem (M) 
and (P) are in one to one correspondence. 
Another natural questions that arises is if $D$ 
inherits natural symmetries that $\Omega$ possesses, and if given $\Omega$, 
does there exist a unique optimal configuration $D$. 
The answer is negative on general domains unless $\Omega$
 has very strong topological restrictions. Symmetrization and 
rearrangement invariant integral methods allow one to prove the next theorem.
\begin{th}{\sc ([CGIKO])}
\label{th2}
 Assume $\Omega$ is symmetric and convex with respect to the hyperplane
$\{x_1=0\}$.
That is for each 
fixed  $x'=(x_2, \cdots, x_n)$ the set 
$$
\{ x_1: (x_1, x')\in\Omega\}
$$
is either empty or an interval of the form $(-c, c)$. 
 Then any optimal 
solution $(u, D)$ is symmetric with respect to the hyperplane
$\{x_1=0\}$ and 
$u$ is decreasing in $x_1$ for $x_1 \ge 0$.
\end{th}
Theorem~\ref{th2} implies the next corollary, which is the only uniqueness 
result proved in [CGIKO].
\begin{cor}
\label{cor3}
Let $\Omega=\{ |x| <1\}$ be the ball. Then the optimal configuration 
is unique for any $\alpha, A$ and furthermore $D$ is an annular region, 
$$
D=\{ x; r(A) < |x| <1\}.
$$
\end{cor}
On other domains we encounter the phenomenon of 
symmetry breaking. Specifically in [CGIKO] we show symmetry breaking 
phenomena on annular domains in ${\bf R}^2$ and 
on dumbbell shaped domains. We have 
\begin{th}{\sc ([CGIKO])}
\label{th4}
Fix any $\alpha >0$ and $\delta\in (0,1)$. Let 
$\Omega_a=\{x\in {\bf R}^2; a < |x| < a+1\}$. 
Then there exists $a_0=a_0(\alpha, \delta)$, such that whenever 
$a > a_0$ and $D$ is an optimal configuration for $\Omega_a$ with 
parameters $\alpha$ and $A=\delta|\Omega_a|$, then $D$ is not rotationally 
symmetric.
\end{th}
Because $\Omega_a$ is rotationally invariant, Theorem~\ref{th4}
 implies that there 
are infinitely many choices for the optimal configuration $D$ 
on annuli. On dumbbell shaped domains we have, in addition to 
symmetry breaking, some extra information on $D$. We define the dumbbell 
shaped domain $\Omega_h$ by 
\begin{equation}
\label{eq6}
\Omega_h=B_1((-2,0))\cup B_1((2,0))\cup 
((-2,2)\times (-h, h)),
\end{equation}
where $B_r(p)=\{x\in {\bf R}^2; |x-p| <r\}$. 
We call the disks $B_r(p)$, the lobes of the dumbbell and the strip 
$(-2,2)\times (-h, h)$ the handle. We have
\begin{th}{\sc ([CGIKO])}
\label{th5}
For any given $\alpha >0$ and $A\in (0, 2\pi)$, there exists $h_0>0$  
such that for domains $\Omega_h$ of (\ref{eq6}) with $h< h_0$, 

\noindent
(a)\ Any optimal pair $(u,D)$ is not symmetric with respect to the 
$x_2-$axis.

\noindent
(b)\ If $A > \pi$, then for any optimal pair $(u,D)$, $D^c$ is totally 
contained in one of the lobes $B_1((\pm2, 0))$. 
\end{th}
We end our summary of results from [CGIKO] with a theorem on convex domains. 
\begin{th}{\sc ([CGIKO])}
\label{th6} 
Suppose $\Omega$ is convex and has a smooth boundary. Then there exists 
$\alpha_0(\Omega, A) >0$, such that 
for any $\alpha < \alpha_0$ and any optimal configuration $D$, one has 

\noindent
(a)\ $\partial D\cap \Omega$ is real-analytic. 

\noindent
(b)\ $D^c$ is convex.
\end{th}
Theorem~\ref{th6}
 should be compared with the basic theorem of Brascamp-Lieb [BL]  
that establishes the convexity of level lines of the first eigenfunction on 
convex domains, when $\alpha=0$. 
Theorem~\ref{th6} extends
the result of [BL] to some values of $\alpha >0$, but it is 
completely open if Theorem~\ref{th6} extends to the case of all $\alpha >0$.

%% file: cm2.tex

We now turn to the results that we prove in this paper.
Our focus will primarily be on general domains and in particular 
on the free boundary for the optimal pair $(u, D)$. 
In a general domain $\Omega$ the first eigenfunction 
$\psi$ (with standard $L^2$ normalization) that is, 
\begin{eqnarray}
\label{eq7}
-\Delta \psi &=&\mu_1\psi\quad \mbox{on}\quad \Omega \\
\psi&=&0\quad \mbox{on}\quad \partial\Omega, 
\quad \int_{\Omega}\psi^2=1,\nonumber
\end{eqnarray}
is real-analytic in the interior and this places very strong 
restrictions on the exceptional sets, i.e. places on the level sets of 
$\psi$, where $\nabla \psi=0$.
If we were in ${\bf R}^2$, the exceptional set would consist of points. 
In this analysis unique continuation plays a role,
since we easily see that if $w=\psi_{x_i}$ then from (\ref{eq7}), $-\Delta
w=\mu_1w$, and thus unique continuation yields some information
on the zero set of $w$. 
 Thus in our problem 
it is clear the free boundary $\{ u=t\}$ on a general domain will 
possess an exceptional set, since $\psi$ in general has one, but 
any attempt to understand the fine structure of the exceptional set, 
Hausdorff measure, rectifiability  etc., 
through a unique continuation approach  
is difficult. The reason being, unique continuation will not 
apply, since $u$ is only weakly regular. This prevents
us from obtaining an equation for $w_1=u_{x_i}$. In addition a further
difficulty is that (\ref{eq5}) is an
equation of the type $-\Delta u+V(x)u=0$, 
and this prevents us from obtaining a homogeneous equation 
that is satisfied by $w_2=u-t$ to which we may apply unique
continuation to
study the level surface  $\{ u=t\}, t\neq 0$. 

Another approach is to view our problem (\ref{eq5}) as a perturbation 
in $\alpha$ from the problem (\ref{eq7}). This approach suffers from 
the fact that we do not get additional information for large $\alpha$. 
One may view this again as a difficulty arising from lack of 
continuation properties in $\alpha$. Two results in this direction are 
proved in [CGIKO]. We reproduce the statements and proof here.
\begin{th}
\label{th7}
For $s\ge 0$, let $[\Omega]^s=\{\psi\le s\}$, where 
$\psi$ is the normalized first eigenfunction of problem (\ref{eq7}). 
Fix $A\in [0, |\Omega|]$ and choose $t_{\Omega}$ such that 
$|[\Omega]^{t_{\Omega}}|=A$. Then for any $\delta >0$, there is 
$\alpha_0=\alpha_0(\delta, \Omega)$ such that if 
$\alpha < \alpha_0$ and $D$ is an optimal configuration for $(\alpha, A)$, 
then $|t-t_{\Omega}|< \delta$ and 
$$
[\Omega]^{t_{\Omega}-\delta}\subset D \subset 
  [\Omega]^{t_{\Omega}+\delta}.
$$
\end{th}
The basic lemma that is used to prove Theorem~\ref{th7}
 can be used to analyse the 
limiting behavior as $A \to |\Omega|$. We have 
\begin{th}
\label{th8}
Let $\Omega$ be a smooth bounded domain. 
Let $\alpha >0$ be fixed. Let $\psi$ be the function of (\ref{eq7}) and 
let 
$$
M=\max_{\Omega}\psi.
$$
Then for any $\delta >0$, there is $A_0=A_0(\delta, \alpha,
\Omega)<|\Omega|$, such 
that whenever $A > A_0$ and $D$ is an optimal configuration for 
$(\alpha, A)$, then 
$$
D^c\subset \{ \psi > M-\delta\}.
$$
\end{th}
The meaning of Theorem~\ref{th7} is that the free boundary for our 
optimization problem, that is the set $\{ u=t\}$, 
is ``trapped'' between the levels $t_{\Omega}-\delta$ and 
$t_{\Omega}+\delta$ of the first eigenfunction $\psi$ 
for the domain $\Omega$. 
However, this information is too weak to conclude anything fine about the 
free boundary even for small $\alpha >0$. 
Theorem~\ref{th8} on the other hand indicates that 
as $A\to |\Omega|$, $D^c$ coalesces onto the set where $\psi$ 
achieves its maximum, $\psi$ being the first eigenfunction on $\Omega$, 
see (\ref{eq7}). Now, keeping in mind part (b) of Theorem~\ref{th5}, 
it is likely that $D^c$ may coalesce onto a strict subset of the 
set where $\psi$ achieves its maximum. 

Even though we have been unable to apply unique continuation 
to study the free boundary for large $\alpha$, it is still possible to apply 
the Hopf lemma [GT, Lemma~3.4] and get some information on the free 
boundary.

A typical result we prove is:
\begin{th}
\label{th9}
Let $\alpha\ge \Lambda(\alpha, A)$. Let ${\cal F}=\{ u=t\}$  
denote the free boundary set. Then there is a subset ${\cal E}$ of 
${\cal F}$ such that 

\noindent
(a) \ ${\cal E}$ is a $G_{\delta}$ set.

\noindent
(b)\ ${\cal F}\setminus {\cal E}$ is a real-analytic, 
$n-1$ dimensional sub-manifold of ${\bf R}^n$.

\noindent
(c) \ If moreover $\alpha > \Lambda(\alpha, A)$, then for every 
$x_0\in {\cal F}$ and every $\epsilon >0$, the ball $B_{\epsilon}(x_0)$ 
contains points of both 
$\{ u >t\}$ and $ \{ u < t\}$. 
\end{th}
We refer to the set ${\cal E}$ as the exceptional set. From the
construction of the exceptional set, we will deduce 
further geometric information regarding the free boundary. 
This is the content of Proposition~\ref{prop9'}, 
which we do not state here (see  section~2). 

We now give a sufficient condition that ensures that the hypothesis 
of Theorem~\ref{th9}, $\alpha\ge \Lambda(\alpha, A)$ is fulfilled.
\begin{prop}
\label{prop10}
Let $\alpha > \mu_1(\Omega)$, where $\mu_1(\Omega)$ is the first 
Dirichlet eigenvalue for $-\Delta$ on $\Omega$. 
Then there exists $A_0=A_0(\alpha)$, such that $\alpha\ge \Lambda(\alpha, A)$ 
for all $A < A_0$.
 Furthermore, for $C_1=||\psi||_\infty^{-2}$ 
and for fixed $A\in (0, C_1)$, there exists $\alpha_0$  
such that  $\alpha\ge \Lambda(\alpha, A)$ 
for all $\alpha \ge \alpha_0$.
\end{prop}
We lastly investigate the effect of the curvature of $\partial\Omega$ 
on the free boundary and the ``thickness'' of the optimal configuration. 
As observed earlier in the remarks after Theorem~\ref{th1},
 $D$ always contains 
a tubular neighborhood of the boundary. The theorem that follows 
demonstrates in the model case of an annulus in ${\bf R}^2$, that 
at places where $\partial\Omega$ has large ``negative'' curvature, 
one finds that $D$ is ``thin''. To state our result we will set 
up some notation. 
Let
$$
\Omega_{\epsilon}=\{ x\in {\bf R}^2; \epsilon < |x| < 1\}, 
\quad B=\{ x\in {\bf R}^2; |x| < 1\}.
$$
For fixed $\epsilon_0$, let $A < \pi (1-\epsilon_0^2).$ 
For $\epsilon < \epsilon_0$ and any fixed $\alpha >0$, 
 let $(u_{\epsilon}, D_{\epsilon})$ denote the optimal pair for 
$\Omega_{\epsilon}$, with lowest eigenvalue 
$\Lambda_{\Omega_\epsilon}(\alpha,A)$, with 
constraint $|D_{\epsilon}|=A$. We now claim that for $\mu_1(B)$, the
first Dirichlet eigenvalue of the unit disk, we have
\begin{equation}
\label{eq8}
\Lambda_{\Omega_{\epsilon}}(\alpha, A)>
\mu_1(B).
\end{equation}
Since $\Omega_{\epsilon} \subset B$, using  $(u_{\epsilon}, D_{\epsilon})$ 
as a trial pair in the variational characterization for 
$\Lambda_{\Omega_\epsilon}(\alpha,A)$, we have
$$
 \Lambda_{\Omega_{\epsilon}}(\alpha, A)
=\int_{\Omega_{\epsilon}}|\nabla u_{\epsilon}|^2 
+\alpha\int_{\Omega_{\epsilon}}\chi_{D_{\epsilon}}u_{\epsilon}^2
> \int_{\Omega_{\epsilon}}|\nabla u_{\epsilon}|^2\geq \mu_1(B).
$$
This establishes, (\ref{eq8}). 
As a consequence of (\ref{eq8}), imposing the hypothesis 
$\alpha \leq \mu_1(B)$, ensures that 
for every $\epsilon\ge 0$, 
$\alpha < \overline{\alpha}_{\Omega_{\epsilon}}(A)$ (see the definition 
(\ref{eq4}) ). Next, if $\alpha < \overline{\alpha}_{\Omega_{\epsilon}}(A)$, 
and if   $D_{\epsilon}$ is radially distributed,
Theorem~2 of [CGIKO], yields that $D_{\epsilon}$ has the form, 
\begin{equation}
\label{eq8'}
D_{\epsilon}=\{ x\in {\bf R}^2; \epsilon < |x| < r_{\epsilon}\ \ {\rm
or}\  
\quad R_{\epsilon} < |x| < 1\}
\end{equation}
for some $r_{\epsilon}, R_{\epsilon}$, 
$\epsilon <r_{\epsilon} < R_{\epsilon} < 1$. Thus if $\alpha \leq
\mu_1(B)$ we 
may assume that if 
$D_{\epsilon}$ is radially distributed, 
then the set $D_{\epsilon}$ has the form 
described by (\ref{eq8'}) for every $\epsilon>0$. We have
\begin{th}
\label{th11}
Assume  $\alpha \leq \mu_1(B)$, and $A>0$ is prescribed. 
Given this choice of $\alpha$ and $A$, let 
$(u_{\epsilon}, D_{\epsilon})$ be an optimal pair with $|D_{\epsilon}|=A$. 
Assume $D_{\epsilon}$ is radially distributed and hence of the form 
(\ref{eq8'}). Then,
$$
\lim\sup_{\epsilon\to 0}r_{\epsilon}=0.
$$
\end{th}
Thus the implication  is that $D_{\epsilon}$ thins out on the boundary layer 
in contact with $\{ x; |x|=\epsilon\}$, the inner boundary of 
$\partial\Omega_{\epsilon}$. As $\epsilon\to 0$, the curvature 
of $\{ x; |x|=\epsilon\}$ is increasing and negative as seen from 
$\Omega_{\epsilon}$. Thus in this model case one may conclude that 
$\mbox{diam}\ (D)$ is small on parts of $D$ which are in contact with 
pieces of $\partial\Omega$, where the curvature of $\partial\Omega$ is 
large and where $\partial\Omega$ is concave as seen from $\Omega$. 

The paper [CGIKO] discusses the historical antecedents of this problem and the interested reader is referred to this paper for a discussion. 
Furthermore the optimization problem discussed here, is amenable 
to being modelled on a computer. Details of the numerical simulation are 
available in [CGIKO] and the interested reader may find the source of the 
algorithms used and the shape of the optimal configuration in many 
types of domains obtained by these numerical studies.

%% file: cm3.tex
\section{Proofs of the Theorems}
In this section we prove Theorems~\ref{th7}-\ref{th11} and 
Proposition~\ref{prop10}. 
We begin with the proof of Theorem~\ref{th7}. 
We need a preparatory Lemma, that is well-known in perturbation theory and in 
the Physics literature [B, Appendix 39, p. 469]. 
We want a slightly more precise form, though the technique of proof is 
standard and the basic idea follows from [B].
\begin{lem}
\label{lem1}
Fix $D\subset \Omega$. Let $u_{\alpha, D}$ be the ( positive, 
$L^2$ normalized ) first eigenfunction of $-\Delta +\alpha\chi_D$ 
with eigenvalue $\lambda(\alpha, D)$. Then there is a constant 
$C=C_{\Omega}$ such that for $0\le \alpha\le 1$ 
($\psi, \mu_1$ refers to (\ref{eq7})), 
\par
\noindent
(a)\quad $0\le \lambda(\alpha, D)-\mu_1\le \alpha,$
\par
\noindent
(b)\ \quad $\| u_{\alpha, D}-\psi\|_{H^2(\Omega)}\le C\alpha,$
\par
\noindent
(c) \quad  $\| u_{\alpha, D}-\psi\|_{L^{\infty}(\Omega)}\le C\alpha.$
\end{lem}
\noindent
{\bf Proof of Lemma 1:}\ 
Recall  we have set $u_{\alpha, D}=u$. 
Note 
\begin{eqnarray*}
\lambda(\alpha, D)\le \int_{\Omega}(|\nabla \psi|^2 +\alpha\chi_D\psi^2)
&\le &\mu_1 +\alpha\int_{\Omega}\psi^2\\
&\le& \mu_1 + \alpha.
\end{eqnarray*}
Thus, $\lambda(\alpha, D) -\mu_1\le \alpha.$ Next, 
$$
\lambda(\alpha, D)=
\int_{\Omega}(|\nabla u|^2 +\alpha\chi_D u^2)
\ge \mu_1 +\alpha \int_{\Omega}\chi_D u^2.
$$
Thus, $\lambda(\alpha, D) -\mu_1\ge 0$, and we have (a). 
To prove (b), let $\{\psi_k\}_{k=1}^{\infty}$ be an orthogonal basis 
of eigenfunction of $-\Delta$ with Dirichlet boundary conditions (Note 
$\psi_1=\psi$ of problem (\ref{eq6})). 
The corresponding eigenvalues will be denoted by $\{ \mu_k\}_{k=1}^{\infty}$, 
where it is well-known that $\mu_1$ is simple. Expanding $u$, we have 
$u=\sum_{j=1}^{\infty}\beta_j\psi_j$, and thus 
$
(-\Delta -\mu_1)u=\sum_{j=2}^{\infty}\beta_j(\mu_j-\mu_1)\psi_j
$ 
and  
\begin{equation}
\label{eq9}
\|(\Delta +\mu_1) u\|_2^2 =\sum_{j=2}^{\infty}\beta_j^2(\mu_j-\mu_1)^2,
\end{equation}
where $\|\cdot\|_2$ denotes the $L^2(\Omega)$ norm. 
From $-\Delta u +\alpha\chi_D u=\lambda(\alpha, D) u$, we get 
$$
(-\Delta -\mu_1)u= (\lambda(\alpha, D) -\mu_1)u-\alpha\chi_D u.
$$
Therefore applying (a), 
$\|(\Delta +\mu_1) u\|_2\le C\alpha$. Since $\mu_1$ is simple, 
there exists $\delta >0, \delta=\delta(\Omega)$, such that 
$\mu_j-\mu_1 \ge \delta >0$ for $j\ge 2$. 
Then from (\ref{eq9}) we get 
\begin{equation}
\label{eq10}
\sum_{j=2}^{\infty}\beta_j^2\le \frac{C\alpha^2}{\delta^2}.
 \end{equation}
We re-write $u$ as $u=\beta_1\psi_1 +\Psi$, and (\ref{eq10}) gives 
$\|\Psi\|_2\le C\alpha\delta^{-1}$. Now  
$1=\|u\|_2 =\beta_1^2+\|\Psi\|_2^2$, thus
$$
|\beta_1-1|\le \frac{\|\Psi\|_2^2}{1+\beta_1}\le \frac{C\alpha^2}{\delta^2}.
$$
Here we used the fact that $\beta_1 =(u,\psi_1) >0$, because both $u$ and 
$\psi_1$ are positive. Therefore, 
\begin{equation}
\label{eq11}
\| u-\psi_1\|_2^2=(\beta_1-1)^2+\|\Psi\|_2^2\le 
 \frac{C\alpha^2}{\delta^2}\le C\alpha^2.
\end{equation}
All the remaining consequences follow from (\ref{eq11}). 
From (\ref{eq6}), 
\begin{equation}
\label{eq12}
-\Delta (u-\psi_1)=(\lambda(\alpha, D)-\mu_1)u 
+\mu_1 (u-\psi_1)-\alpha\chi_D u.
\end{equation}
We re-write (\ref{eq12}) as 
$$
-\Delta(u-\psi_1)-\mu_1 (u-\psi_1)=
(\lambda(\alpha,D)-\mu_1)u -\alpha\chi_D u=g.
$$
Now from [GT, Theorem~8.15] again, it follows that 
$$
\| u-\psi_1\|_{\infty}\le C\| u-\psi_1\|_2 + C\alpha.
$$
Using (\ref{eq11}) on the right side, 
$$
\| u-\psi_1\|_{\infty}\le  C\alpha,
$$
which is part (c). Using $\|g\|_{\infty}\le C\alpha$  and part (c), we 
conclude $\|\Delta (u-\psi_1)\|_{\infty}\le C\alpha$, which is (b). 
$\Box$

Theorem~\ref{th7} and \ref{th8} are now consequences of Lemma~\ref{lem1}.

\noindent
{\bf Proof of Theorem~\ref{th7}:}\ 
Apply Lemma~\ref{lem1} (c) to the optimal pair $(u, D)$. Choose $\alpha_0=
\delta/(2C)$, so that 
$\|u-\psi_1\|_{\infty}\le \delta/2$ for 
$\alpha \le \alpha_0$. From Theorem~\ref{th1}, if $x\in D, u(x) \le t$, and 
so $\psi(x) \le t+\delta/2$ and hence $D\subset [\Omega]^{t+\delta/2}$. 
In a similar way we establish $[\Omega]^{t-\delta/2}\subset D$. 
Thus we have 
$$
[\Omega]^{t-\delta/2}\subset D \subset [\Omega]^{t+\delta/2}.
$$
From the statement above, we get 
$|[\Omega]^{t-\delta/2}|\le A\le |[\Omega]^{t+\delta/2}|$, 
and thus by continuity, there exists $t_{\Omega}$ such that 
$A=|[\Omega]^{t_{\Omega}}|$, and $|t_{\Omega}-t| < \delta/2$. 
From this assertion the assertions of Theorem~\ref{th7} follow. 
$\Box$

\noindent
{\bf Proof of Theorem~\ref{th8}:}\ 
We begin by showing that a slight modification of the 
proof of Lemma~\ref{lem1} yields, 
\begin{equation}
\label{eq13}
\|u-\psi_1\|_{\infty}\le C_{\alpha, \Omega}(|\Omega|-A).
\end{equation}
We show first, 
\begin{equation}
\label{eq14}
|\mu_1- (\Lambda(\alpha, A)-\alpha)|\le  C_{\alpha, \Omega}(|\Omega|-A).
\end{equation}
We re-write our equation for $u$ as
\begin{equation}
\label{eq15}
-\Delta u -\alpha\chi_{D^c}u=(\Lambda-\alpha)u, \quad \Lambda=
\Lambda(\alpha, A).
\end{equation}
From (\ref{eq15}) we have
$$
\mu_1 -\alpha\int_{\Omega}\chi_{D^c}u^2\le \int_{\Omega}|\nabla u|^2 
-\alpha\int_{\Omega}\chi_{D^c}u^2=\Lambda -\alpha.
$$
Thus, 
$$
\mu_1-(\Lambda -\alpha)\le \alpha\int_{\Omega}\chi_{D^c}u^2
\le C|D^c|=C(|\Omega|-A).
$$
Next, we have
$$
\Lambda -\alpha \le \int_{\Omega}(|\nabla \psi|^2-\alpha \chi_{D^c}\psi^2)
\le \mu_1 -\alpha\int_{\Omega}\chi_{D^c}\psi^2
$$
which yields $0\le \mu_1 -(\Lambda -\alpha)$. 
The assertion (\ref{eq14}) follows. Using (\ref{eq14}) and 
the equation (\ref{eq15}) we can proceed as in Lemma~\ref{lem1} to obtain 
(\ref{eq13}). 
If $|\Omega|-A < \delta/(2C_{\alpha, \Omega})$, from (\ref{eq13}) we see 
$[\Omega]^{t-\delta/2}\subset D=\{ u\le t\}\subset [\Omega]^{t+\delta/2}$. 
Since   $|D|=A$, we have 
$A \le |[\Omega]^{t+\delta/2}|=|\{ \psi\le t+\delta/2\}|$. 
Thus if in addition $A > A_0$, we can arrange the situation so as to have 
$M-t\le \delta/2$. 
So $[\Omega]^{M-\delta}\subset [\Omega]^{t-\delta/2}\subset D$. 
The conclusion  
$[\Omega]^{M-\delta}\subset D$ is readily seen to be equivalent 
to the assertion made in Theorem~\ref{th8}. 
$\Box$

We need some preparatory lemmas before we prove Theorem~\ref{th9}. 
As usual $u=u_{\alpha, D}$ will denote the solution to our 
optimization problem and consequently $u$ will satisfy (\ref{eq5}). 
\begin{lem}
\label{lem2}
(a)\ Fix any $\alpha >0$. Let the free boundary set be ${\cal F}$, 
${\cal F}=\{ u=t\}$. We let $D^+=\{ x; u(x) > t\}$. 
Assume $D^+$ satisfies an interior sphere condition with respect to 
$x_0\in {\cal F}$, that is there exists a ball $B$, 
$B\subset D^+$ and $\partial B\cap {\cal F}=\{x_0\}$. 
Then $|\nabla u(x_0)|\neq 0$. 

\noindent
(b)\ Let $D^-=\{x; u(x) <t\}$ (in the situation 
$\alpha\neq \overline{\alpha}_{\Omega}(A)$, $D^-=D$). 
Assume that $\alpha\ge \Lambda_{\Omega}(\alpha, A)$. Let $D^-$ satisfy 
an interior  sphere condition with respect to 
$x_0\in {\cal F}$, that is there exists a ball $B$, 
$B\subset D^-$ and $\partial B\cap {\cal F}=\{x_0\}$. 
Then $|\nabla u(x_0)|\neq 0$. 
\end{lem}
\noindent
{\bf Proof of Lemma~\ref{lem2}:}\ 
The proof of both parts of our lemma rely on Hopf's lemma [GT, Lemma~3.4]. 
We prove (a). Set $\phi=t-u$. We observe that in the ball 
$B\subset D^+$, 
$u$ satisfies from (\ref{eq5}) 
$$
-\Delta u=\Lambda(\alpha, A) u.
$$
Thus $\Delta \phi =\Lambda(\alpha, A) u\ge 0$ in $B$, and 
$\phi< 0$ on $B$ with $\phi(x_0)=0$. Hopf's lemma then yields 
$|\nabla \phi(x_0)|=|\nabla u(x_0)|\neq 0$. 

(b). \ The proof of this part is similar to part (a). Since $B \subset D^-$, 
from (\ref{eq5}) we see on $B$ we have 
$$
-\Delta u + \alpha u=\Lambda u.
$$
Since $\alpha\ge \Lambda$, we easily see $\Delta u\ge 0$ on $B$. 
Thus on $B\subset D^-$ we have $\Delta \phi\le 0, \phi >0$ on $B$ and 
$\phi(x_0)=0$. Thus Hopf's lemma again yields 
$|\nabla \phi(x_0)|=|\nabla u(x_0)|\neq 0$. 
$\Box$

\begin{lem}
\label{lem3}
Let $h(\eta, p), \eta\in {\bf R}, p=(p_1, \cdots, p_n)\in {\bf R}^n$ be a locally 
bounded function. Let $w\in C^1(\Omega)$ satisfy 
\begin{equation}
\label{eq16}
\Delta w =h(w, \nabla w).
\end{equation}
Assume furthermore $h$ is smooth in the variable $p$. Assume at the point 
$x_0\in\Omega$, $\nabla w(x_0)\neq 0$. Then there exists a ball $B$, 
$x_0\in B$, such that the set, 
$$
\{ x\in B; w(x)=w(x_0)\}={\cal S}
$$
is a smooth hypersurface of ${\bf R}^n$. 
If in addition $h$ is real-analytic in the variable $p$, the set 
${\cal S}$ is also real-analytic.
\end{lem}
\noindent
{\bf Proof of Lemma 3:}\ 
  Since $h(w, \nabla w)$ is locally bounded, it follows by elliptic 
estimates that $w\in C^{1, \gamma}\cap W^{2,s},\ s<\infty$. Thus by the implicit function theorem,
 since $\nabla w(x_0)\neq 0$, we conclude that ${\cal S}$ is a 
$C^{1, \gamma}$ hypersurface for all $\gamma <1$. Now we shall improve 
the regularity of the hypersurface ${\cal S}$. By a rotation of coordinates 
we may assume $w_{x_i}(x_0)= 0, i=1, \cdots, n-1$ and 
$w_{x_n}(x_0)\neq 0$. Let $x'=(x_1, \cdots, x_{n-1})$ and consider the map,
$$
\Psi:\  B \to {\bf R}^n, \quad \Psi(x', x_n)=(x', w(x', x_n)).
$$
We denote points in the image of $\Psi$, by 
$y=(y', y_n)$ where $y'=(y_1,\cdots, y_{n-1})$ and $y_n=w(x',x_n)$. 
Let $\Psi^{-1}$ denote the inverse map to $\Psi$, which will exist if $B$ 
is picked to be small. We have,
 $$
\Psi^{-1}(y', y_n)=(y', F(y', y_n)).
$$
Now, 
$$
F(x', w(x', x_n))=x_n.
$$
Differentiating the equation above we get the equations,
\begin{equation}
\label{eq17}
F_{y_i}+w_{x_i}F_{y_n}=0, i=1, \cdots, n-1, \quad \mbox{and} \quad 
F_{y_n}w_{x_n}=1.
\end{equation}
By the chain rule, 
$$
\frac{\partial}{\partial x_i}=
\frac{\partial}{\partial y_i}+w_{x_i}\frac{\partial}{\partial y_n}, 
\quad 
\frac{\partial}{\partial x_n}=
w_{x_n}\frac{\partial}{\partial y_n}.
$$
From (\ref{eq17}), 
$w_{x_i}=-F_{y_i}/F_{y_n}, i=1, \cdots, n-1$ and 
$w_{x_n}=1/F_{y_n}$. Thus, 
\begin{eqnarray*}
\Delta w&=& 
\sum_{i=1}^{n-1}
\biggl( \frac{\partial}{\partial y_i}+w_{x_i}\frac{\partial}{\partial y_n} 
\biggr)
\biggl( \frac{-F_{y_i}}{F_{y_n}}\biggr)
+w_{x_n}\frac{\partial}{\partial y_n}\biggl(
\frac{1}{F_{y_n}}\biggr)\\
&=& 
\sum_{i=1}^{n-1}
\biggl( \frac{\partial}{\partial y_i} -
\frac{F_{y_i}}{F_{y_n}}\frac{\partial}{\partial y_n} 
\biggr)
\biggl( \frac{-F_{y_i}}{F_{y_n}}\biggr)
+\frac{1}{F_{y_n}}\frac{\partial}{\partial y_n}\biggl(
\frac{1}{F_{y_n}}\biggr)\\
     &=&LF.
\end{eqnarray*}
Next we freeze the coefficients of $L$ at $x_0$. 
Since 
$w_{x_i}(x_0)=0, i=1, \cdots, n-1$ and $w_{x_n}(x_0)=a\neq 0$, we see
from (\ref{eq17}), $F_{y_i}=0$ and $F_{y_n}=a^{-1}$ at $\Psi(x_0)$. 
At $y_0=\Psi(x_0)$, we have 
\begin{eqnarray*}
LF&=& -\frac{1}{F_{y_n}}\sum_{i=1}^{n-1}F_{y_iy_i}
-\frac{1}{F_{y_n}^3}F_{y_ny_n}\\
&=& 
 -\frac{1}{F_{y_n}} \biggl[
\sum_{i=1}^{n-1}F_{y_iy_i}
+\frac{1}{F_{y_n}^2}F_{y_ny_n}\biggr].
\end{eqnarray*}
Thus if $B$ is picked suitably small, $LF$ is an elliptic, 
quasi-linear operator. Since $w$ satisfies (\ref{eq16}), $F$ satisfies 
\begin{equation}
\label{eq18}
LF=h\biggl(y_n, {\cal A}(\frac{-F_{y_i}}{F_{y_n}}, \frac{1}{F_{y_n}})\biggr)
\end{equation}
where ${\cal A}$ is a fixed matrix in $O(n)$, the rotation group, 
associated with our rotation of coordinates, and 
${\cal A}(z_i, z_n)$, denotes the product of ${\cal A}$ and the vector 
$(z_1, \cdots, z_n), (z_i)=(z_1, \cdots, z_{n-1}).$ 
Since $w\in C^{1,\gamma}\cap W^{2,s}$ for all $\gamma <1, s<\infty$, 
we see the coefficients of 
$L$ are in $C^{0,\gamma}\cap W^{1,s},\ s<\infty$. 
Differentiating equation (\ref{eq18}) 
in any of the variables $y_i, i\neq n$, we may apply a standard
bootstrap argument using well-known elliptic 
estimates for example [GT, Theorem~9.11] to conclude from the fact that 
$h$ is smooth in the variable $p$, that $F$ is smooth in any of 
the variables  $y_i, i=1,\cdots, n-1$. If one has in addition that $h$ 
is real-analytic in the variables $p$, then applying the results of 
Morrey [M, Theorem~C] or Friedman [F, Theorems~1, 4], we can conclude 
that $F$ is real-analytic in the variables $y_i, i=1,\cdots, n-1$. 
Now the defining equation for ${\cal S}$ is given by 
$x_n=F(x_1, \cdots, x_{n-1}, w(x_0))$. Since $F$ is smooth (real-analytic) 
depending on the regularity of $h$ in the $p$ variables, it follows 
that ${\cal S}$ is a smooth (real-analytic) manifold, depending on the 
fact that  $h$ is smooth (real-analytic) in the $p$ variables. 
$\Box$

We are now ready to prove Theorem~\ref{th9}.

\noindent
{\bf Proof of Theorem~\ref{th9}:}\ 
We construct the exceptional set ${\cal E}$. Let 
$$
K_n=\{x\in \Omega;\  \mbox{distance}\ (x, {\cal F})=\frac{1}{n}\}.
$$
Now define 
$$
F_n=\{ x\in {\cal F}; \    
\mbox{distance}\ (x, K_n)=\frac{1}{n}\}.
$$
The sets $K_n$ and $F_n$ are closed sets for all $n\in {\bf N}$. We define the 
exceptional set by 
$$
{\cal E}={\cal F}\setminus ( \cup_{n=1}^{\infty} F_n).
$$
Since each set ${\cal F}\cap F_n^c$ is open in ${\cal F}$, 
${\cal E}$ is a $G_{\delta}$ set. This proves Theorem~\ref{th9} (a). 

We shall now prove that for each point $x_0\in {\cal F}\setminus {\cal E}$, 
we can construct either an interior ball $B \subset D^+$, 
such that $\partial B \cap {\cal F}=\{x_0\}$ or 
an interior ball $B \subset D^-$, 
such that $\partial B \cap {\cal F}=\{x_0\}$. 
Since we are assuming $\alpha\ge \Lambda(\alpha, A)$, Lemma~\ref{lem2} 
ensures that for each $x_0\in {\cal F}\setminus {\cal E}$, 
$\nabla u(x_0)\neq 0$. Now the PDE satisfied by $u$, that is (\ref{eq5}), 
can be written as 
$$
\Delta u =h(u),
$$
where $h(\eta)=-\Lambda(\alpha, A)\eta +\alpha\chi_G(\eta)\eta, 
G=\{ \eta; \eta\le t\}$. 
Thus $h\in L^{\infty}_{loc}({\bf R})$ and the hypotheses of 
Lemma~\ref{lem3} apply to $u$. Thus there exists a ball $B_0$, 
centered at $x_0$, such that, ${\cal F}\cap B_0$ is a 
real analytic manifold. So we now verify our claims regarding the interior 
spheres. Fix a point $x_0\in {\cal F}\setminus {\cal E}$. 
Then $x_0\in F_n$ for some $n$. Let $z_0\in K_n$, such that 
$|x_0-z_0|= 1/n$. We claim the ball $B_{1/n}(z_0)$ is totally 
contained in $D^+$ or $D^-$. 
Suppose there are points $z_1\in D^+$ and $z_2\in D^-$ in 
$\overline{B_{1/n}(z_0)}$. 
Then by the continuity of $u$, the line segment joining $z_1$ to $z_2$ 
which also lies inside $B_{1/n}(z_0)$ will contain a point of 
${\cal F}$. Thus distance $(z_0, {\cal F})< 1/n$ and hence 
$z_0\not\in K_n$. Now pick a ball $B', B'\subset B_{1/n}(z_0)$, and 
$B'$ centered along the radius joining $x_0$ to $z_0$ and with 
$\partial B'\cap \partial B_{1/n}(z_0)=\{x_0\}$. 
Clearly $\overline{B'\setminus\{x_0\}}$ is contained either 
in $D^+$ or $D^-$ and $\overline{\partial B'}\cap {\cal F}=\{x_0\}$. 
The hypotheses of Lemma~\ref{lem2} are fulfilled. Theorem~\ref{th9} (b) 
now follows. 

To prove Theorem~\ref{th9} (c), we can argue via the strong maximum 
principle. Since $\alpha > \Lambda(\alpha, A)$, 
applying Theorem~\ref{th1} (d), we see that 
$D=\{ u <t\}$. If for some $\epsilon >0$, the ball $B_{\epsilon}(x_0)$, 
$x_0\in {\cal F}$ contains no points of $D^-$ and only points 
$\{u\ge t\}$, then on $B_{\epsilon}(x_0)$, $u$ satisfies 
$$
-\Delta u =\Lambda(\alpha, A) u\ge 0.
$$  
So $\Delta u\le 0$ on $B_{\epsilon}(x_0)$ and $u\ge t$ on 
$B_{\epsilon}(x_0)$ with $u(x_0)=t$ which is a 
contradiction. We may argue as in Lemma~\ref{lem2}, part (b) and show on 
$B_{\epsilon}(x_0)$ there are also points of $D^+$ 
for every $\epsilon >0$. 
$\Box$

To discuss further geometric properties of ${\cal E}$, we introduce,
$$
K_{\epsilon}=\{x\in \Omega;\  \mbox{distance}\ (x, {\cal F})=\epsilon\}, 
$$
$$
F_{\epsilon}=\{ x\in {\cal F}; \    
\mbox{distance}\ (x, K_{\epsilon})=\epsilon\}.
$$
We have 
\begin{prop}
\label{prop9'}
For every $\alpha >0$, 

\noindent
(a)\ $F_{\epsilon_1}\subset F_{\epsilon_2}$ for $\epsilon_1 > \epsilon_2$.

As a consequence of (a), we have 

 \noindent
(b)\ $\cup_{n\ge 1}F_n =\cup_{\epsilon >0}F_{\epsilon}$, 
$F_n\subset F_m$ for $ m >n$, and 
${\cal E}={\cal F}\setminus \cup_{\epsilon >0}F_{\epsilon}$. 

\noindent
(c) \ If $\alpha > \Lambda(\alpha, A)$, then for every point 
$x_0\in {\cal E}$, any ball $B_{\epsilon}(x_0)$ contains points 
$y_+\in D^+$ and $y_-\in D^-$ such that 
distance\  $(y_{\pm}, {\cal F}) < |y_{\pm}-x_0|$. 

\noindent
(d) \ Furthermore, if $z_0\in {\cal F}$, such that $|y_+-z_0|=$\
distance\ $(y_+,{\cal F})$, then $D^+$ satisfies an interior sphere
condition in the sense of Lemma~\ref{lem2} with respect to $z_0$. Thus by the
proof of Theorem~\ref{th9}, there is a neighborhood $B$ of $z_0$,
such that $\{x\in B: u(x)=u(z_0)\}$ is a real-analytic hypersurface. 
A similar statement holds for $y_-$.
\end{prop}
The meaning of (b) is that since 
${\cal E}={\cal F}\cap (\cap_{n=1}^{\infty}F_n^c)$ and 
$F_n^c \supset F_m^c$ for $m>n$, the exceptional set is really a 
$G_{\delta}$ set formed by the intersection of the nested open sets $F_n^c$. 
The meaning of (c), (d) is that the behaviour of the free boundary in 
the neighborhood of the exceptional set at least in ${\bf R}^2$ is
that of isolated singularities with a conical structure, with
the cone having its vertex at 
$x_0\in {\cal E}$. The cone locally divides ${\bf R}^2$ into at least
two components, one component is contained 
in $D^+$ and the other is contained in $D^-$. This geometric picture is only
heuristic since it still needs to be rigorously established that the
component of ${\cal F}$ that contains $z_0$ also contains $x_0\in {\cal
E}$. Only then can we conclude there is a true conical singularity at $x_0$.

\medskip
\noindent
{\bf Proof of Proposition~\ref{prop9'}:}\ 
Fix $x\in F_{\epsilon_1}$. Then by definition, 
one can find $z\in K_{\epsilon_1}$, such that 
$|z-x| = \epsilon_1$. Since $z\in K_{\epsilon_1}$, the ball 
$B_{\epsilon_1}(z)$ will contain no points of ${\cal F}$ in the interior. 
Now consider the radial line connecting $z$ and $x$, and locate on 
this line a point $y$ such that $|y-x|=\epsilon_2$. 
The ball $B_{\epsilon_2}(y)\subset B_{\epsilon_1}(z)$, and $x$ is the 
sole point in ${\cal F}$ on $\partial B_{\epsilon_2}(y)$. Now by definition 
$y\in K_{\epsilon_2}$, and since $x\in {\cal F}$ and 
$|x-y|=\epsilon_2$, $x\in F_{\epsilon_2}$. 
This proves (a), and (b) is then an elementary consequence of set theory. 

Next, by Theorem~\ref{th9} (c) we know that the ball $B_{\epsilon}(x_0)$ 
with $x_0\in {\cal E}$ ( in general for any point in ${\cal F}$ actually), 
contains points $y_{\pm}$ in $D^{\pm}$ respectively. 
Let $|y_+-x_0|=\delta$. 
We claim distance\ $(y_+, {\cal F}) < \delta$. 
If distance\ $(y_+, {\cal F})\ge \delta$, then $y_+\in K_{\tau}, 
\tau=$ distance\ $ (y_+, {\cal F}), \tau\ge \delta.$  
Since $|y_+-x_0|=\delta\le \tau$, it follows that 
$x_o\in F_{\tau}$. Thus $x_0\not\in {\cal E}$. This proves part (c), 
since an argument similar to the one above takes care of 
the point $y_-\in D^-$.

Lastly we prove part (d). Now, $|y_+-z_0|=\tau$, $z_0\in {\cal F}$. 
Thus by the argument employed in Theorem~\ref{th9}(b) it is easily
seen, that the open ball $B_\tau(y_+)$ contains only points of $D^+$.
Again employing the argument of Theorem~\ref{th9} (b), we can find a
ball $B^\prime\subset B_\tau(y_+)$ such that $\partial B^\prime\cap{\cal
F}=\{z_0\}$ and hence $B^\prime$ is the desired interior ball.
$\Box$

%% file: cm4.tex

We now prove Proposition~\ref{prop10}.

\noindent
{\bf Proof of Proposition~\ref{prop10}:}\ 
Arguing as in Lemma~\ref{lem1}(a), but now using the fact that 
$\|\psi\|_{\infty}\le C$, [GT, Theorem~8.15], 
we have 
$$
\Lambda(\alpha, A)\le \int_{\Omega}|\nabla \psi|^2+\alpha 
\int_{\Omega}\chi_D\psi^2 
\le \mu_1 +\alpha\|\psi\|_{\infty}^2A.
$$
Thus,
\begin{equation}
\label{eq19}
\Lambda(\alpha, A)\le \mu_1 +\alpha C_0 A, \quad 
C_0=\|\psi\|_{\infty}^2.
\end{equation}
Since $\alpha > \mu_1$, we can find $A_0 >0$ such that 
$\alpha C_0A_0\le \alpha - \mu_1$. 
It follows from (\ref{eq19}) that for $A < A_0$, 
$\alpha \ge \Lambda(\alpha, A)$. 
The second part of Proposition~\ref{prop10} also follows from (\ref{eq19}). 
Select $C_1=C_0^{-1}$. If $A < C_1$, $C_0A=1-\epsilon, \epsilon>0$. 
Thus for $\alpha > \alpha_0$, $\mu_1\le \epsilon\alpha$, and hence by 
(\ref{eq19}), $\Lambda(\alpha, A) < \alpha$.
$\Box$

\noindent
{\bf Proof of Theorem~\ref{th11}:}\ 
We set $\Lambda_{\epsilon}=\Lambda_{\Omega_{\epsilon}}(\alpha, A)$, 
$\Lambda=\Lambda_B(\alpha, A)$. $D$ will denote the optimal configuration 
for $B$, that is $D=\{ x; r_0 < |x| < 1\}, \pi(1-r_0^2)=A$. 
We also need to consider the first Dirichlet eigenvalue $\lambda_{\epsilon}$, 
of the problem,
\begin{eqnarray*}
-\Delta g +\alpha\chi_D g&=& \lambda_{\epsilon}g \quad \mbox{in}\quad 
\Omega_{\epsilon}\\
g|_{\partial\Omega_{\epsilon}}=0.
\end{eqnarray*}
We claim 
\begin{equation}
\label{eq20}
0\le \Lambda_{\epsilon}-\Lambda \le C|\log\epsilon|^{-1}.
\end{equation}
Extending $u_{\epsilon}$ to $\{ |x| < \epsilon\}$ by setting 
$u_{\epsilon}=0$ on $\{ |x| < \epsilon\}$ and using this extended 
function as a trial function with $D_{\epsilon}$ as a trial configuration 
on $B=\{ |x| <1\}$, we see 
$\Lambda\le \Lambda_{\epsilon}$ and so we have the left side in 
(\ref{eq20}). 
Next by Theorem~2 in Swanson [S], 
\begin{equation}
\label{eq21}
0\le \lambda_{\epsilon}-\Lambda \le c |\log\epsilon|^{-1}.
\end{equation}
In fact $0\le \lambda_{\epsilon}-\Lambda$ is  simply a consequence of 
domain monotonicity. By the variational characterization  
$\Lambda_{\epsilon}\le \lambda_{\epsilon}$, and thus from (\ref{eq21})  
we easily have (\ref{eq20}). 

We now establish $\lim\sup_{\epsilon\to 0}r_{\epsilon}=0$, by 
contradiction. 
Assume there is a sequence $\epsilon_j \searrow 0$, 
$\lim_{j\to\infty}r_{\epsilon_j}=\delta >0$. 
Then $R_{\epsilon_j}\to
b_{\delta}$ and since $|D_{\epsilon_j}|=A$,
the limit set ${\cal D}_{\delta}$ is:
$$
{\cal D}_{\delta}=\{ x; |x|< \delta\ \ {\rm or}\quad b_{\delta} < |x| <1\}, 
\quad 
|{\cal D}_{\delta}|=A.
$$
We set,
$$
D_{\delta,\epsilon}=\{ x; \epsilon < |x|< \delta\ \ {\rm or}\quad b_{\delta} < |x| <1\}.
$$
We use the notation $\lambda_{\epsilon}(D_{\delta, \epsilon})$ for the 
first Dirichlet eigenvalue on $\Omega_{\epsilon}$ for the problem,
\begin{eqnarray}
\label{eq22}
-\Delta w &+&\alpha\chi_{D_{\delta,\epsilon}}w
=\lambda_{\epsilon}(D_{\delta,\epsilon}) w \quad \mbox{on}\quad 
\Omega_{\epsilon}\\
w&=&0 
\quad \mbox{on}\quad 
\partial\Omega_{\epsilon}, \quad \int_{\Omega_{\epsilon}}w^2=1.
\nonumber
\end{eqnarray}
$\lambda({\cal D}_{\delta})$ will denote the first Dirichlet eigenvalue for 
$-\Delta + \alpha\chi_{{\cal D}_{\delta}}$ on $B$. We claim,
\begin{equation}
\label{eq22'}
|\Lambda_{\epsilon_j} -\lambda({\cal D}_{\delta})|\to 0\quad \mbox{as}\quad 
\epsilon_j\to 0.
\end{equation}
We have
\begin{eqnarray*}
|\Lambda_{\epsilon_j} -\lambda({\cal D}_{\delta})|
&\le&
|\Lambda_{\epsilon_j} -\lambda_{\epsilon_j}(D_{\delta,\epsilon_j})|+
|\lambda_{\epsilon_j}(D_{\delta,\epsilon_j}) -\lambda({\cal D}_{\delta})|\\
&=& J_1+J_2.
\end{eqnarray*}
By using Theorem~2 in [S], one can conclude 
$J_2\le C|\log \epsilon_j|^{-1}\to 0$ as $\epsilon_j\to 0$. 
We now show $J_1\to 0$. 
Let $(u_{\epsilon}, D_{\epsilon})$ be the optimal pair for 
$\Omega_{\epsilon}$. Then, using the eigenfunction $w$ from (\ref{eq22}), 
and the uniform bounds 
$\|w\|_{L^{\infty}(\Omega_{\epsilon})}\le C$, 
$C$ independent of $\epsilon >0$, which follows from [GT, Theorem~8.15], 
we have 
\begin{eqnarray*}
\Lambda_{\epsilon}&\le& \int_{\Omega_{\epsilon}}
(|\nabla w|^2 +\alpha\chi_{D_{\epsilon}}w^2)\\
&=&
\int_{\Omega_{\epsilon}}
(|\nabla w|^2 +\alpha\chi_{D_{\delta, \epsilon}}w^2)
+\alpha \int_{\Omega_{\epsilon}}(\chi_{D_{\epsilon}}
-\chi_{D_{\delta, \epsilon}} )w^2\\
&\le& \lambda_{\epsilon}(D_{\delta,\epsilon})+C\alpha
|D_{\epsilon}\bigtriangleup  D_{\delta,\epsilon}|, 
\end{eqnarray*}
where $D_{\epsilon}\bigtriangleup  D_{\delta,\epsilon}$ is the symmetric 
difference of the sets 
$D_{\epsilon}, D_{\delta,\epsilon}$. Thus for the sequence $\epsilon_j$, 
we easily have  
$|D_{\epsilon}\bigtriangleup  D_{\delta,\epsilon}|\to 0$ as $j\to\infty$. 
We conclude
$$
\Lambda_{\epsilon_j} \le 
\lambda_{\epsilon_j}(D_{\delta,\epsilon_j}) +o(1), \quad j\to\infty.
$$
Likewise using the function $u_{\epsilon}$ in the argument above, we have 
\begin{eqnarray*}
\lambda_{\epsilon_j}(D_{\delta,\epsilon_j})&\le&
 \Lambda_{\epsilon_j} +C\alpha
|D_{\epsilon_j}\bigtriangleup  D_{\delta,\epsilon_j}|\\
&\le& \Lambda_{\epsilon_j}+ o(1), \quad j\to \infty.
\end{eqnarray*}
Thus, 
$$
|\lambda_{\epsilon_j}(D_{\delta,\epsilon_j})-
 \Lambda_{\epsilon_j}|\to 0 \quad\mbox{as}\quad j\to\infty.
$$
Thus $J_1\to 0$ as $\epsilon_j\to 0$. 
This establishes (\ref{eq22'}). We infer from (\ref{eq20}) and 
(\ref{eq22'}) that as $j\to\infty$, 
$$
|\Lambda -\lambda({\cal D}_{\delta})|\le 
|\Lambda-\Lambda_{\epsilon_j}|+
|\Lambda_{\epsilon_j}-\lambda({\cal D}_{\delta})|\to 0.
$$
Hence $\Lambda=\lambda({\cal D}_{\delta})$. However, this contradicts 
Corollary~\ref{cor3} of our introduction, a proof of which is supplied in 
[CGIKO]. $\Box$

%% file: cm5.tex
\section{Open Problems and Conjectures}
A number of open problems and conjectures can be stated based on the numerical 
data in [CGIKO] and the rigorous results there and also on results 
proved here. We will outline some.
\begin{prob}: (Uniqueness of the optimal configuration)\ 
The only domain for which we have established the uniqueness of the optimal configuration is the ball, see Corollary~\ref{cor3}. 
Is $D$ unique if $\Omega$ is convex?
\end{prob}
\begin{prob}:(Continuation Problem)\ 
Theorem~\ref{th6} states that on a convex domain $\Omega$, $D^c$ is convex 
for small $\alpha>0$. Is it possible to continue along the values $\alpha >0$,
 to obtain convexity of $D^c$ for all $\alpha>0$?
\end{prob}
\begin{prob}:(The free boundary on general domains)\ 
The free boundary ${\cal F}$ on general domains will contain 
an exceptional set ${\cal E}$ as constructed in the proof of
Theorem~8. 
What is the Hausdorff dimension of ${\cal E}$? 
Is at least ${\cal F}$ a rectifiable set?
One suspects ${\cal E}$ consists of points, where real-analytic arcs intersect 
if $\Omega\subset {\bf R}^2$. 
\end{prob}
\begin{prob}:(Monotonicity of $D$)\ 
Suppose $A < A'$, then does this imply $D_{\alpha, A} \subset 
D_{\alpha, A'}$? If symmetry breaking occurs this statement needs to 
be modified. Nevertheless on domains where the 
optimal configuration is unique, 
does the above monotonicity hold?
\end{prob}
\begin{prob}:(Symmetry breaking on annuli)\ 
When $\Omega$ is an annulus, what is the shape of $D$ precisely?
The proof of Theorem~\ref{th4} in our introduction as presented in 
[CGIKO] and the numerical computations presented in [CGIKO] 
suggest that $D^c$ lies between two rays $\theta=0$ and $\theta=\beta$. 
In fact the results in [CGIKO] suggest that, 
$\beta=\pi/N, N=N(\alpha, |D|/|\Omega|)$, and 
$N\to \infty$ as $|D|\to |\Omega|$. 
\end{prob}
\begin{prob}:(Influence of the boundary curvature)\ 
In Theorem~\ref{th11} we saw in a model case that the diameter of $D$ is 
affected by the curvature of $\partial\Omega$. Investigate this phenomena 
on general domains $\Omega$.   
\end{prob}
We refer the interested reader to [CGIKO] for further problems and conjectures.

%% file: cm-ref.tex

\begin{center}
{\bf REFERENCES}
\end{center}

\vspace{5mm}

[B] M. Born, Atomic Physics, Dover Publications, 1989. 

[BL] H. J. Brascamp, E. Lieb, On extensions of the Brunn-Minkowski
and Prekopa-Leindler theorems, including inequalities for log concave 
functions, and with an application to the diffusion equation,
J. Funct. Analysis, 22(1976), 366-389.

[CGIKO] S. Chanillo, D. Grieser, M. Imai, K. Kurata and I. Ohnishi, 
Symmetry breaking and other phenomena in the optimization of eigenvalues of
composite membranes, to appear.

[F] A. Friedman, On the regularity of the solutions of non-linear elliptic and 
parabolic systems of partial differential equations,
 Journ. Math. Mech., 7(1958), 43--60.

[GT] D. Gilbarg and N. Trudinger, 
Elliptic Partial Differential Equations of Second Order, Springer-Verlag, 1983.

[M] C.B. Morrey, On the analyticity of solutions of analytic 
non-linear elliptic systems of PDE I, 
Amer. J. of Math. 80(1958), 198--218.

[S] C.A. Swanson, Asymptotic variational formulae for eigenvalues, Canad. Math. Bull., Vol.6, No.1(1963), 15--25.

\bigskip
\bigskip
\bigskip